\documentclass[a4paper,11pt]{article}
\usepackage[english]{babel}
\usepackage[a4paper]{geometry}
\usepackage{hyperref}
\usepackage{bm}
\usepackage{authblk}
\usepackage{siunitx}
\usepackage{graphicx}
\usepackage{makeidx}        
\usepackage{subcaption}  
\usepackage{graphicx}                                    
\usepackage{multicol}        
\usepackage[bottom]{footmisc}
\usepackage{hyperref}
\usepackage{lineno}
\usepackage{algorithm}
\usepackage{algorithmic}
\usepackage{amsmath} 
\usepackage{amssymb}
\usepackage[square,numbers]{natbib}

\title{\Large Model Reduction for Transport-Dominated Problems via Cross-Correlation Based Snapshot Registration}

\author[1]{Harshith~Gowrachari\footnote{hgowrach@sissa.it}}
\author[2]{Giovanni~Stabile\footnote{giovanni.stabile@santannapisa.it}}
\author[1]{Gianluigi~Rozza\footnote{grozza@sissa.it}}

\affil[1]{\small Mathematics Area, mathLab, International School for Advanced Studies,
via Bonomea 265, 34136 Trieste, Italy}
\affil[2]{\small Biorobotics Institute, Sant’Anna School of Advanced Studies, V.le R. Piaggio 34, 56025, Pontedera, Pisa, Italy.}

\date{}

\begin{document}

\maketitle
\begin{abstract}

 Traditional linear approximation methods, such as proper orthogonal decomposition and the reduced basis method, are ill-suited for transport-dominated problems due to the slow decay of the Kolmogorov \textit{n}-width, leading to inefficient and inaccurate reduced-order models. In this work, we propose a model reduction approach for transport-dominated problems by employing cross-correlation based snapshot registration to accelerate the Kolmogorov \textit{n}-width decay, thereby enabling the construction of efficient and accurate reduced-order models using linear approximation methods. We propose a complete framework comprising offline-online stages for the development of reduced order models using the cross-correlation based snapshots registration. The effectiveness of the proposed approach is demonstrated using two test cases: 1D travelling waves and the higher-order methods benchmark test case, 2D isentropic convective vortex. 
 \end{abstract}


\section{Introduction}
Reduced order modelling is a rapidly advancing field in computational science and engineering. Industries have a keen interest in obtaining reduced order models (ROMs) for engineering systems, particularly in applications involving control, optimization, and uncertainty quantification. ROMs offer significant computational efficiency in many query scenarios and are suitable for real-time computations. These models achieve computational efficiency by approximating high-dimensional parametric partial differential equations (PDEs) solutions with low-dimensional representations, substantially reducing the computational cost compared to full-order models (FOMs). \\

To obtain the FOM solutions of the parametric PDEs describing real-world scenarios, discretization methods such as the finite element method (FEM) and finite volume method (FVM) are widely employed. Although these methods ensure accuracy, they are computationally expensive. ROM framework mitigates this challenge through a two-stage approach:
\begin{itemize}
    \item Offline (training) stage -- This computationally expensive phase involves generating full-order model (FOM) solutions (snapshots) and constructing the ROM by approximating the high-dimensional solution manifold with a low-dimensional approximation subspace.
    \vspace{2mm}
    \item Online (testing and predictive) stage -- This significantly more efficient phase utilizes the low-dimensional approximation subspace obtained during the offline stage to rapidly predict reduced solutions for new, unseen parameters, substantially reducing computational expenses.
\end{itemize}

Several model reduction techniques are available to obtain the low-dimensional approximation subspaces \cite{quarteroni2015reduced, hesthaven2016certified, Chinesta2017, rozza2022advanced}. Among these, linear approximation methods such as proper orthogonal decomposition (POD) and reduced basis (RB) methods are widely used and highly effective in various applications. However, certain problems of interest, such as transport-dominated problems, exhibit slow decay of Kolmogorov \textit{n}-width (KnW). This slow decay of KnW limits the model reduction achievable through linear subspace approximations \cite{Pinkus1985}. \\

The solution manifolds of the transport-dominated problems, such as travelling waves, convective vortex, and vortex shedding, exhibit slow decay of the Kolmogorov \textit{n}-width (KnW), making it difficult to obtain a low-dimensional linear approximation subspace for the construction of accurate and efficient ROMs. The KnW, as defined in \cite{Pinkus1985}, provides a rigorous measure of the reducibility of the solution manifold \(\mathcal{M}\) using linear approximation subspaces. \\

For a solution manifold \(\mathcal{M}\) comprising \textit{f} elements (solutions across all time instances) embedded in a normed linear space \( (X_{\mathcal{N}}, \| \cdot \|_{X_{\mathcal{N}}}) \), the Kolmogorov \textit{n}-width \(d_n(\mathcal{M})\) is defined as:

\begin{equation}
    d_n\left(\mathcal{M}\right) := \inf _{E_n \subset X_N} \sup _{f \in \mathcal{M}} \inf _{g \in E_n}\|f-g\|_{X_{\mathcal{N}}}
\end{equation}

where \(E_n\) is a linear subspace of \(X_{\mathcal{N}}\) with dimension \(n\). The KnW \(d_n(\mathcal{M})\) represents the maximum approximation error obtained when any element \(f \in \mathcal{M}\) is approximated by an element \(g \in E_n\). Furthermore, as faster the \(d_n(\mathcal{M})\) decay with increasing \(n\), indicates the effectiveness of an \textit{n} dimensional low-dimensional linear subspace is in approximating the solution manifold  \(\mathcal{M}\). \\

Recent advancements in model-reduction techniques have focused on addressing the challenges posed by slow Kolmogorov \textit{n}-width decay, being the main reason for the development of non-linear approaches \cite{peherstorfer2022breaking}. These methods include performing model reduction on non-linear manifolds, such as constructing non-linear trial subspaces using convolutional autoencoders (CAEs) \cite{kashima2016nonlinear, hartman2017deep, crisovan2019model, hoang2022projection, lee2020model, fresca2021comprehensive, kim2022fast, romor2023non}, and employing quadratic approximation manifolds \cite{Jain2017, barnett2022quadratic, Geelen2023}. Other strategies involve incorporating adaptive enrichment techniques,  \cite{haasdonk2008adaptive, washabaugh2012nonlinear, peherstorfer2015online, peherstorfer2020model, bruna2024neural}, and transforming linear subspaces through transport maps using registration methods \cite{ohlberger2013nonlinear, Mojgani2017LagrangianBM, iollo2014advection, mojgani2017arbitrary, cagniart2017model, rim2018transport, taddei2020registration, torlo2020model, nonino2023reduced, Nonino2024}. These registration based approaches yield a low-dimensional transformed linear subspace that effectively approximates full-order solutions. \\

In this work, we employ cross-correlation based snapshot registration for model reduction of transport-dominated problems. This approach leverages cross-correlation to align multiple snapshots to a reference, accelerating the decay of KnW and yielding a low-dimensional approximation subspace for the construction of efficient ROMs. It is particularly effective in scenarios where conventional linear approximation methods struggle, providing a more robust framework for capturing the system's essential dynamics, and enabling the construction of efficient and accurate ROMs using the transformed linear approximation subspace obtained through registration. \\

In reduced-order modelling, once the reduced approximation subspace is obtained, it is essential to compute the evolution of the dynamics within this subspace. This requires determining the modal coefficients or latent coordinates, which can be achieved through either intrusive or non-intrusive methods. In the intrusive approach, the discretized governing equations (PDEs) are projected onto the low-dimensional subspace, resulting in a system of low-dimensional ODEs. This procedure, known as Galerkin or Petrov-Galerkin projection based ROM, is discussed in \cite{stabile2018finite, stabile2019reduced}. In this work, we focus specifically on non-intrusive methods, where, after obtaining the set of snapshots (FOM solutions), the ROM is constructed using proper orthogonal decomposition (POD) with regression (POD-R) strategy. This method is discussed in detail in the following section. \\

This article is organized as follows, in section \ref{ROM} we discuss the construction of a non-intrusive or data-driven reduced order model using POD with regression (POD-R) strategy; in section \ref{registration} we discuss the employment of the cross-correlations based snapshots registration approach; in section \ref{ROM-via-registration} the implementation of the complete algorithm for the development of ROM utilizing the cross-correlations based registration is mentioned; in section \ref{sec:numerical_experiment} we test our proposed approach by applying to transport-dominated problems, 1D travelling waves and 2D isentropic convective vortex; and in section \ref{sec:Conclusion} the conclusions, perspectives, and future directions are discussed. \\

After working on this manuscript, we came across the work \cite{Issan2023}, which closely aligns with the framework we propose. In their approach, cross-correlation based snapshot registration, a widely used technique in medical imaging \cite{berberidis2002new, chelbi2018features}, is integrated with the operator inference framework \cite{Peherstorfer2016} for the construction of ROMs. This framework is referred to as shift operator inference \cite{Issan2023}, where shifts are obtained via cross-correlation. In this work, we integrate cross-correlation based snapshot registration within the non-intrusive data-driven ROM framework \cite{Tezzele2022}. 

\section{Reduced-order modelling of time-dependent problems}
\label{ROM}
In this work, we are focusing on time-dependent problems, considering the abstract partial differential equation (PDE) as shown given in (\ref{Eqn: abstract_problem}), defined over the domain \(\Omega \subset \mathbb{R}^{n}\) and the time interval \(T = [0, t_f]\):  

\begin{equation}
    \mathcal{L}(s(t)) = 0, \quad t \in T.
    \label{Eqn: abstract_problem}
\end{equation}

Where $s(t)$ is an unknown function, the operator $\mathcal{L}$ incorporates the differential operators, forcing terms and boundary conditions that govern the dynamics of $s(t)$. To obtain the full-order solutions of the given time-dependent problem, here we consider the finite-difference method, but one can choose to use any discretization method, such as the finite-volume method or finite-element method, and later, followed by a time integration scheme. We denote the solutions database as a set of time-snapshot pairs given by \( \{ t_i, s_i \}_{i=1}^{N_T} \) with $t_i \in T \subset \mathbb{R}^{N_T} $ with $N_T$ number of time instances and $ s_i \in \mathbb{R}^{N_{h}}$ with $N_h$ number of degrees of freedom. We represent the snapshots matrix $ S \in \mathbb{R}^{{N_h}\times N_{T}} $ as :

\begin{equation}
S=\left[\begin{array}{cccc}
\mid & \mid & & \mid \\
s_1 & s_2 & \ldots & s_{N_{T}} \\
\mid & \mid & & \mid
\end{array}\right].
\end{equation}

To reduce the computational cost of the associated discretization method, in this section, we present the general framework for non-intrusive ROMs using proper orthogonal decomposition with regression (POD-R) strategy \cite{Tezzele2022}. Here, we first recall the proper orthogonal decomposition, a widely used linear approximation method in the ROM community, which is exploited together with the regression strategy to obtain efficient and accurate ROMs. 

\subsection{Proper orthogonal decomposition}

 We use POD to obtain the reduced linear approximation subspace $U_r$, where $ dim(U_r) = N_r \ll N_h$, which consists of POD modes spanned column-wise. The POD modes are obtained by performing singular value decomposition (SVD) \cite{PODstewart1993early} of the snapshots matrix: 

\begin{equation}
 S = U \Sigma V^T, 
\end{equation}

where, ${U} \in \mathbb{R}^{{N_h}\times N_{h}}$ and ${V} \in \mathbb{R}^{{N_T}\times N_{T}}$ are the two orthonormal matrices. $\mathbf{\Sigma} \in \mathbb{R}^{N_h \times N_T}$ is the diagonal matrix with \textit{r} non-zero singular values, arranged in descending order $\sigma_{1} \geq \sigma_{2} \geq \cdot\cdot\cdot\ \geq \sigma_{r}>0$, which indicates the energy contribution of the corresponding modes. Here, \textit{r} is the rank of the snapshot matrix \textit{S}. The reduced linear approximation subspace $U_r \in \mathbb{R}^{N_h \times N_r} $ is constructed such that :  

\begin{equation}
\mathbf{U}_{\mathbf{r}}=\underset{\hat{\mathbf{U}} \in \mathbb{R}^{N_h \times N_r}}{\operatorname{argmin}} \frac{1}{\sqrt{N_T}}\left\|\mathbf{S}-\hat{\mathbf{U}} \hat{\mathbf{U}}^{T} \mathbf{S}\right\|_F
\end{equation}

where $ \left\|  \cdot \right\|_F $ is a frobenius norm and \textit{S} is the snapshots matrix. \\ 

The POD modes are column entries of the left singular matrix \textit{U} corresponding to the largest singular values. The energy retained by the first \textit{m} modes is given by the ratio of the energy contained by the first $m$ modes and the energy contained by all \textit{r} modes :

\begin{equation}
E(m) = \frac{\sum_{i=1}^m \sigma_i^2}{\sum_{i=1}^r \sigma_i^2}. 
\label{eqn: POD tolerance}
\end{equation}

Generally, the reduced space $U_r$ is constructed by considering the first $N_r$ columns of \textit{U} \cite{eckart1936approximation}, such that the $E(N_r)$ contains a user-defined threshold, usually about $95 - 99.9 \%$ of the cumulative energy. 

\subsection{POD with regression}
In this work, POD with regression (POD-R) strategy is employed to obtain non-intrusive data-driven ROM \cite{Tezzele2022}, which exploits regression to approximate the trajectories of the full-order model. Once we have the reduced space $U_r$, we approximate the solutions of $s(t)$ by reduced expansion: 

\begin{equation}
    s(t) = U_r c_r(t)
\end{equation}

where, $c_r \in \mathbb{R}^{N_r}$ is the temporal expansion coefficients. We compute the temporal expansion coefficients as shown in (\ref{eqn:coeff matrix}), by projecting the snapshot matrix \textit{S} onto the reduced subspace $U_r$.   

\begin{equation}
 c(t_i) = U_r^T \hspace{1mm} s(t_i) \hspace{2mm} \forall \hspace{2mm} i = 1, \dots , N_{train}. 
\label{eqn:coeff matrix}
\end{equation}

Following the previous steps, we can only approximate the snapshots of the initial database (training set) via these temporal expansion coefficients. To predict snapshots for new time instances, we need to build a map from time $t_i$ to $c_i$ temporal expansion coefficients. Here, the main aim is to construct a regression $ \mathcal{I} : T \mapsto \mathbb{R}^{N_r} $ which approximates the map $\mathcal{F}$ as in (\ref{mapping_param_to_coeff}) given by a set of $N_{train}$ input-output pairs $ \left\{t_i, c_i\right\}_{i=1}^{N_{train}}$, where $t_i$ is the time associated with the $i^{th}$ snapshot and $c_i$ are the corresponding temporal coefficients obtained by (\ref{eqn:coeff matrix}).
\begin{equation}
    \mathcal{F} : t \in T \mapsto c \in \mathbb{R}^{N_{r}}. 
    \label{mapping_param_to_coeff}
\end{equation}
To construct this regression $\mathcal{I}$, there are various techniques available, such as linear interpolation, Gaussian process regression (GPR), radial basis function (RBF) interpolation, multi-fidelity methods, and artificial neural networks (ANNs), to name a few. In this work, we exploit radial basis function (RBF) interpolation to devise this map. The constructed regression model is used to predict the snapshot $s^*$ for new time instance $t^*$ via inverse coordinate transformation: \\  
\begin{equation}
    \mathbf{s}(t^*) = U_r \hspace{1mm} I(t^*).
\end{equation}

\section{Cross-correlation based snapshot registration}\label{registration}

Cross-correlation measures the similarity between two signals as a function of the displacement or lag between them. It plays a crucial role in signal processing, pattern recognition, time series analysis, and filtering. Cross-correlation aids in identifying or quantifying the relationships between signals and is used to detect signals that are embedded in noise. Image registration through cross-correlation is an essential and widely adopted technique in the field of medical imaging \cite{berberidis2002new, chelbi2018features}. In this work, we exploit cross-correlation to register the snapshots of transport-dominated problems and aim to accelerate the KnW decay. \newline

\noindent \textbf{Continuous Cross-Correlation:} Let us consider two continuous functions \( f(t), g(t) \in \mathbb{C} \). The continuous cross-correlation between \textit{f} and \textit{g} is defined as:

\begin{equation}
(f \star g)(\tau) \triangleq \int_{-\infty}^{\infty} \overline{f(t)} g(t+\tau) \, dt
\end{equation}

where \( \overline{f(t)} \) denotes the complex conjugate of \( f(t) \), and \( \tau \) is the time lag or displacement. The cross-correlation measures the similarity between \( f(t) \) and a time-shifted version of \( g(t) \), integrating the product of complex conjugate \( \overline{f(t)} \) and the shifted \( g(t + \tau )\) over all time. \\

\noindent \textbf{Discrete Cross-Correlation:} Let us consider two discrete signals \( f[n], g[n] \in \mathbb{C}^{N_h} \), where \( N_h \) is the number of elements in the signals. The discrete cross-correlation between \textit{f} and \textit{g} is defined as:

\begin{equation}
(f \star g)[k] \triangleq \sum_{n=0}^{N_h - 1} \overline{f[n]} g[n+k]
\end{equation}

where \( \overline{f[n]} \) is the complex conjugate of the discrete signal \( f[n] \), \( n \) denotes the index of the signal elements, and \( k \) is the displacement or lag. The discrete cross-correlation is computed by summing the element-wise product between \( \overline{f[n]} \) and the shifted sequence \( g[n+k] \). Depending on the context, \( g[n+k] \) may be assumed to be zero outside its defined range (zero-padding) or periodically extended (circular correlation). \\

Cross-correlation offers a robust framework for analyzing the relationship between time-shifted signals. In this work, we leverage cross-correlation snapshot registration, as outlined in Algorithm \ref{Algorithm1}. In this framework, snapshots are treated as discrete signals, and the cross-correlation between each snapshot $s_i$ and the reference snapshot $s^{\text{ref}}$ is used to compute the shifts required for registration. Here, we employ the registration of snapshots via cross-correlation to circumvent KnW decay for transport-dominated problems. This method aims to determine the optimal shift for each snapshot, aligning them with the reference based on the maximization of their cross-correlation. \\

The process of cross-correlations based registration is as follows: First, we select a set of training snapshots \( \{ s_1, s_2, \dots, s_{N_{train}} \} \) corresponding to time instances \( \{ t_1, t_2, \dots, t_{N_{train}} \} \), and choose one snapshot as the reference \( s^{\text{ref}} \).  For each training snapshot $s_i$, we compute the cross-correlation\footnote{We use \texttt{scipy} python library to compute correlation} between the reference $s^{\text{ref}}$ and the snapshot $s_i$ at various time shifts (lags). Cross-correlation measures the similarity between two signals as a function of their relative time shifts. The optimal shift $\Delta_i$ for each snapshot $s_i$ is determined by identifying the shift that maximizes the cross-correlation, i.e., the shift that yields the highest similarity between the reference and the respective snapshot $s_i$. Once the optimal shift  $\Delta_i$ is determined for each snapshot $s_i$, these shifts are employed to align the snapshots with the reference $s^{\text{ref}}$. Then, we store the registered snapshots $s_i^{\text{ref}}$ and their corresponding optimal shifts for subsequent steps in developing the ROM. 

\begin{algorithm}
\caption{: Registration via cross-correlation}
\begin{algorithmic}[1]
\label{Algorithm1}
\STATE \textbf{Input:} Set of training snapshots \( \{ s_1, s_2, \dots, s_{N_{train}} \} \) for corresponding time instances \( \{ t_1, t_2, \dots, t_{N_{train}} \} \), reference snapshot \( s^{\text{ref}} \)
\STATE \textbf{Output:} Set of Registered snapshots \( \{ s_1^{\text{ref}}, s_2^{\text{ref}}, \dots, s_{N_{train}}^{\text{ref}} \} \), corresponding shifts \( \Delta = \{ \Delta_1, \Delta_2, \dots, \Delta_{N_{train}} \} \)
\vspace{2mm}
\FOR{each snapshot \( s_i \) in set of training snapshots}
\vspace{2mm}
    \STATE Compute discrete cross-correlation between the reference snapshot \( s^{\text{ref}} \) and the snapshot \( s_i \):
    \[
    \text{correlation}(s^{\text{ref}}, s_i) = \text{correlate}(s^{\text{ref}}, s_i)
    \]
    \STATE Obtain optimal shift \( \Delta_i \):
    \[
    \Delta_i = \text{argmax} \left( \text{correlation}(s^{\text{ref}}, s_i) \right) - (N_h - 1)
    \]
    \STATE To register snapshots: Apply the shift \( \Delta_i \) to the snapshot \( s_i \) using circular shift:
    \[
        s_i^{\text{ref}} = \text{shift}(s_i, \Delta_i)
    \]
    \STATE Store the Registered snapshot \( s_i^{\text{ref}} \) and the shift \( \Delta_i \)
\vspace{2mm}
\ENDFOR
\vspace{2mm}
\STATE \textbf{return} the set of Registered snapshots \( S^{\text{ref}} = \{ s_1^{\text{ref}}, s_2^{\text{ref}}, \dots, s_{N_{train}}^{\text{ref}} \} \) and shifts \( \Delta_i \). 
\end{algorithmic}
\end{algorithm}

\section{ROM via cross-correlation based snapshot registration}\label{ROM-via-registration}

In this section, we present the complete workflow for constructing reduced-order models (ROMs) using the cross-correlation based snapshot registration framework. As outlined in Algorithm \ref{Algorithm2}, the workflow consists of two stages: (1) the offline stage, where snapshots are registered and the ROM is constructed, and (2) the online stage, the ROM predictions which are at the reference are reverted to the correct physical frame. \\

In the offline phase, the snapshots are registered using Algorithm \ref{Algorithm1}, which exploits cross-correlation to employ registration and yields a registered snapshots matrix $S^{\text{ref}}$ and a corresponding shifts matrix $\Delta$. Now, all the snapshots are in the reference frame, accelerating the KnW and yielding a low-dimensional linear approximation subspace. We construct ROM\footnote{using \texttt{EZyRB} package \cite{demo18ezyrb}, a Python library for data-driven ROMs.} utilizing the registered snapshots matrix $S^{\text{ref}}$, following the steps mentioned in section \ref{ROM}, with radial basis function (RBF) interpolation used for regression. To obtain shifts for new time instances in the online phase, we construct a map $\alpha$ using a regression on the training pairs $\{ t_i, \Delta_i \}$. In the online phase, the ROM predictions are initially at the reference frame, and we exploit the corresponding optimal shift obtained by the map $\alpha$ to revert the snapshots to the correct physical frame. 

\begin{algorithm}
\caption{: ROM via cross-correlation based snapshot registration}
\begin{algorithmic}[1]
\label{Algorithm2}
\STATE \textbf{Input:} Set of training snapshots \( \{ s_1, s_2, \dots, s_{N_{train}} \} \) for corresponding time instances \( \{ t_1, t_2, \dots, t_{N_{train}} \} \), reference snapshot \( s^{\text{ref}} \), $N_r$ number of POD modes. 
\STATE \textbf{Output:} Prediction of snapshots via cross-correlation based snapshots registration. 
    \vspace{2mm} 
    \underline{\textbf{Offline stage}}
    \vspace{1mm} 
    
    \STATE Registration of snapshots via cross-correlation using Algorithm \ref{Algorithm1}:

    {Obtain the Registered snapshots matrix $S^{\text{ref}}$ and the corresponding optimal shifts matrix $\boldsymbol{\Delta}$ }

    \STATE Construct ROM on the registered snapshots matrix $S^{\text{ref}}$ following the steps in section \ref{ROM}. 

        \[
            \text{rom} = \text{ROM}(S^{\text{ref}}, \text{POD}, \text{RBF})
        \]

    \STATE Construct map $\mathcal{\alpha}$ using regression for $\{ t_i, \Delta_i \}$ pairs: To predict shifts for new time instances in the online phase.

    \vspace{2mm} 
    \noindent \underline{\textbf{Online stage}}
    \vspace{1mm} 

    \FOR{each new time instance \( t_i \) in test set: }
    \STATE $s_i^{\text{ref}}$ = rom.predict($t_i$)
    \STATE Obtain the optimal shift $\Delta_i^{*}$ using map $\alpha$ .  
    \STATE Transport back the predicted snapshots to the correct physical frame by reverting them with the optimal shift:
    \[
        s_i = \text{shift}(s_i^{\text{ref}}, - \Delta_i^{*})
    \]
\ENDFOR
\STATE \textbf{return} the set of predicted snapshots \( \{ s_1 , s_2, \dots, s_{N_{test}} \} \) in correct physical frame. 
\end{algorithmic}
\end{algorithm}

\section{Numerical experiments}
\label{sec:numerical_experiment}
In this section, we evaluate the proposed approach to model reduction using cross-correlation based snapshot registration to construct ROMs for transport-dominated problems. We consider two numerical experiments: the 1D travelling waves and a higher-order benchmark test case - the 2D isentropic convective vortex.

\subsection{1D Travelling wave}

\begin{figure}[h!]
    \includegraphics[scale=0.6]{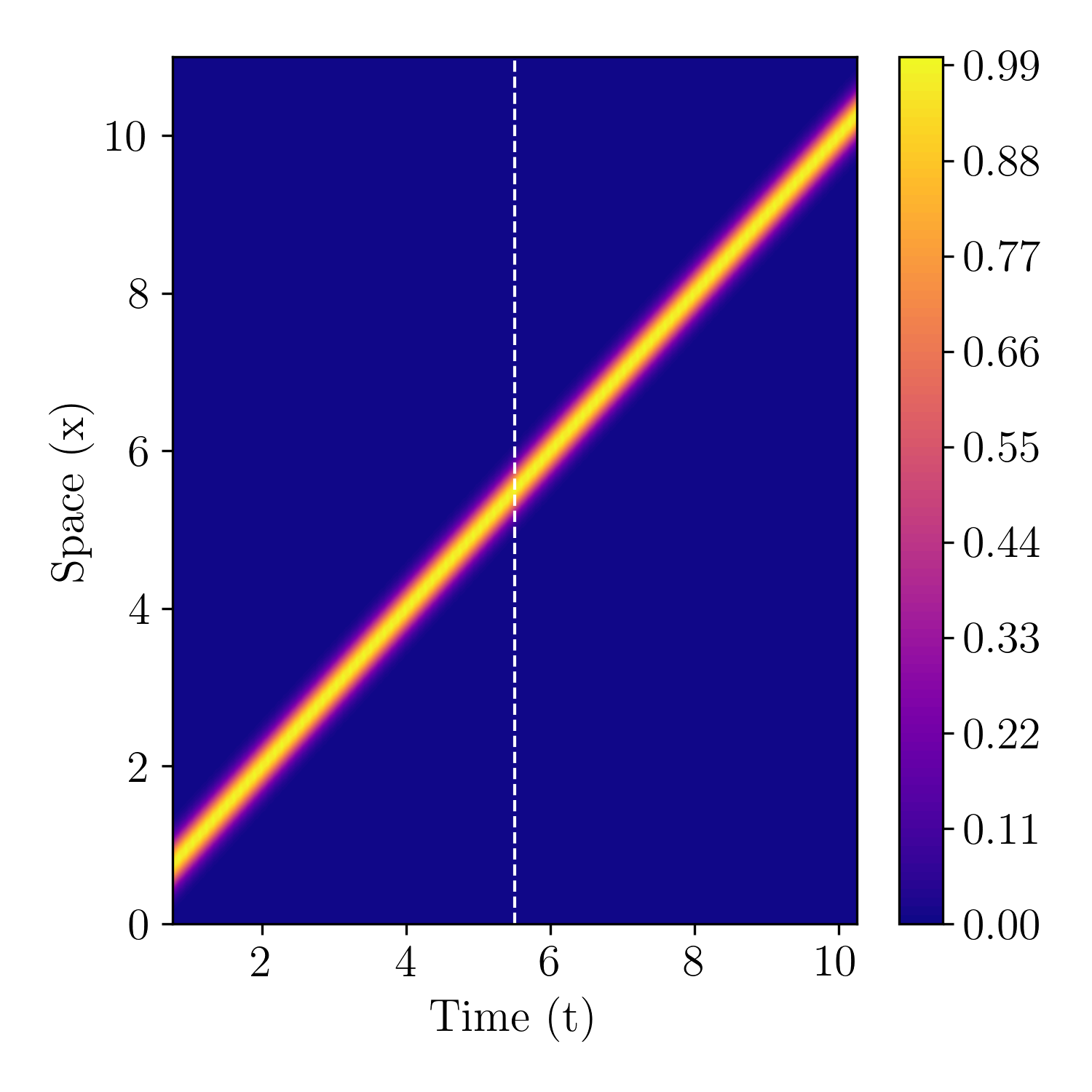}
    \centering
    \caption{2D contour of the original snapshots of the 1D travelling waves. Here, the white line separates the first 50$\%$ dataset as the training set and the rest as the test set.}
    \label{fig:1Doriginal_gaussian_snaps}
\end{figure}

We consider a parameterized 1D Gaussian function (\ref{eq:gaussian}) to mimic the solutions of the linear transport equation considered in \cite{fresca2021comprehensive, Geelen2023, harshithNNsPOD}, with $\beta = 1$, $ x \in [0, 10.25]$ is the spatial coordinate, constant $\sigma$ (standard deviation), and $t \in [0, 10.25]$ is the mean, interpreted here as time to treat this as a parameterized time-dependent problem. We generate a database consisting of 100 time instances and their corresponding snapshots over the 256 equidistant spatial nodes in $x \in [0, 10.25]$. We divide 50$\%$ of the database as a training set and the rest as a testing set, as shown in the Figure (\ref{fig:1Doriginal_gaussian_snaps})

\begin{equation}\label{eq:gaussian}
     f(t)= \beta e^{-(x-{t})^2 /\left(2 \sigma^2\right)}.
\end{equation}

\begin{figure}[h!]
    \centering
    \begin{subfigure}[b]{0.45\textwidth}
        \centering
        \includegraphics[width=\textwidth]{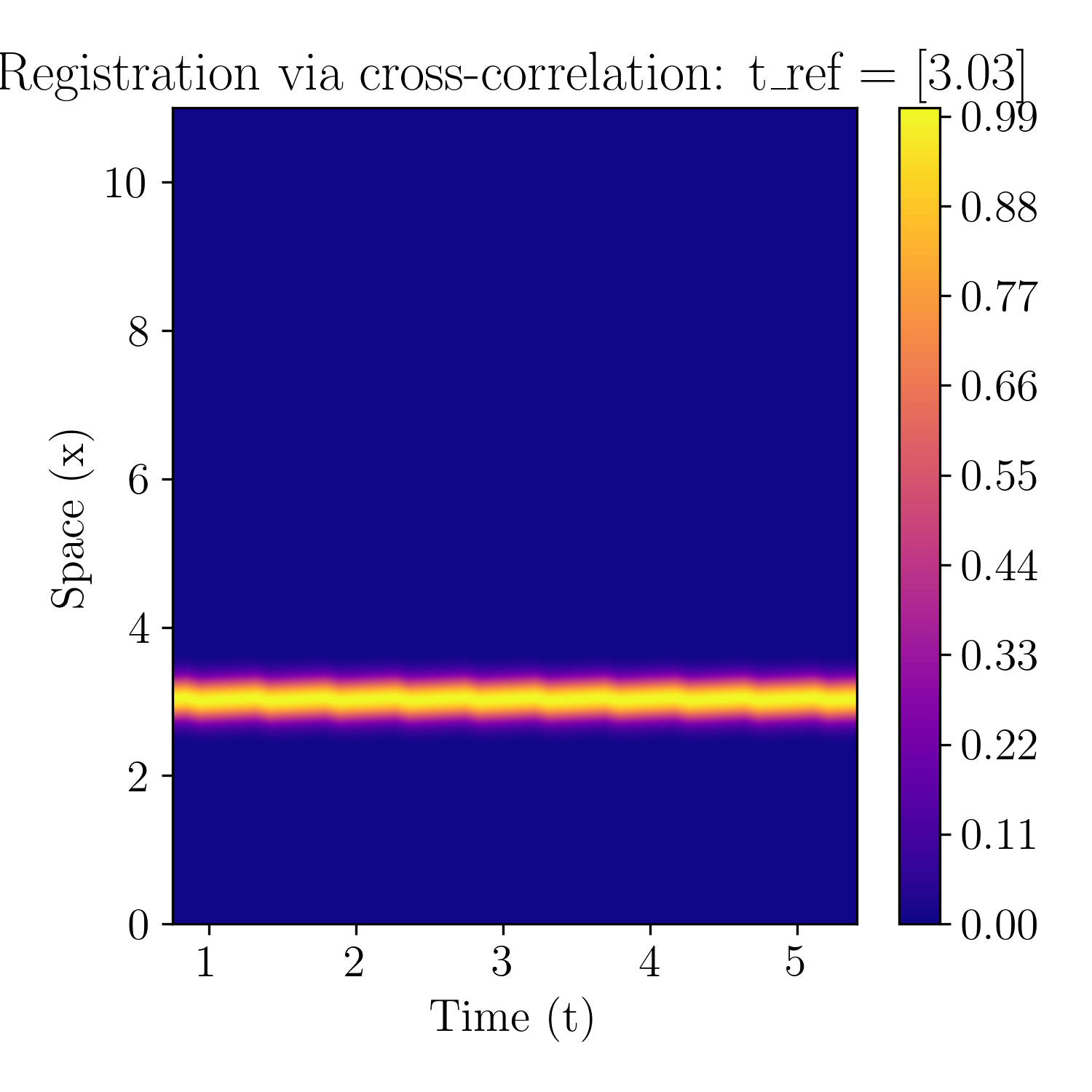}
        \caption{2D contour of the registered train set. }
    \end{subfigure}
    \begin{subfigure}[b]{0.45\textwidth}
        \centering
        \includegraphics[width=\textwidth]{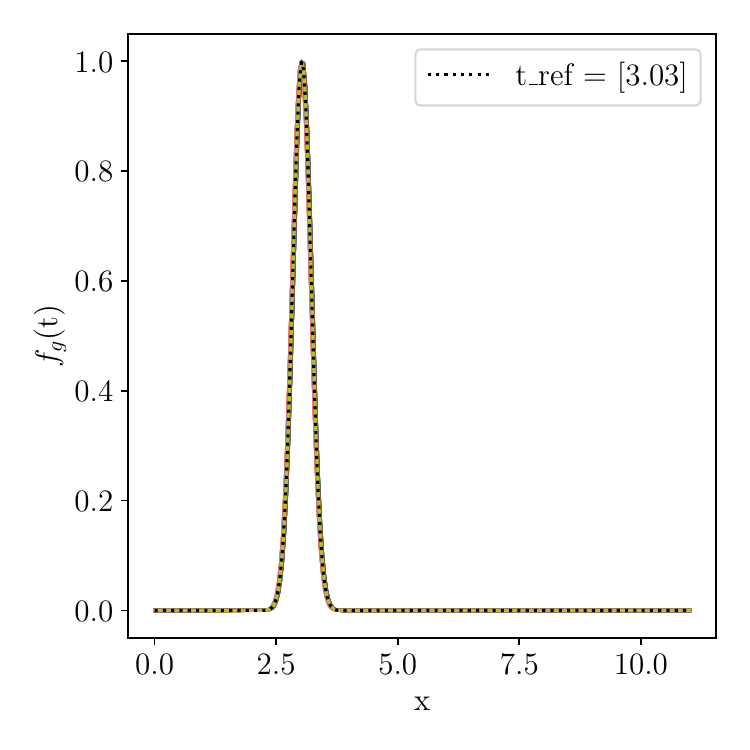}
        \caption{Snapshots aligned with reference.}
    \end{subfigure}
    \caption{1D travelling waves training set registered at the reference configuration. }
    \label{fig:1D_registered_snapshots}
\end{figure}

\begin{figure}[htbp]
    \centering
    \begin{subfigure}[b]{0.48\textwidth}
        \centering
        \includegraphics[width=\textwidth]{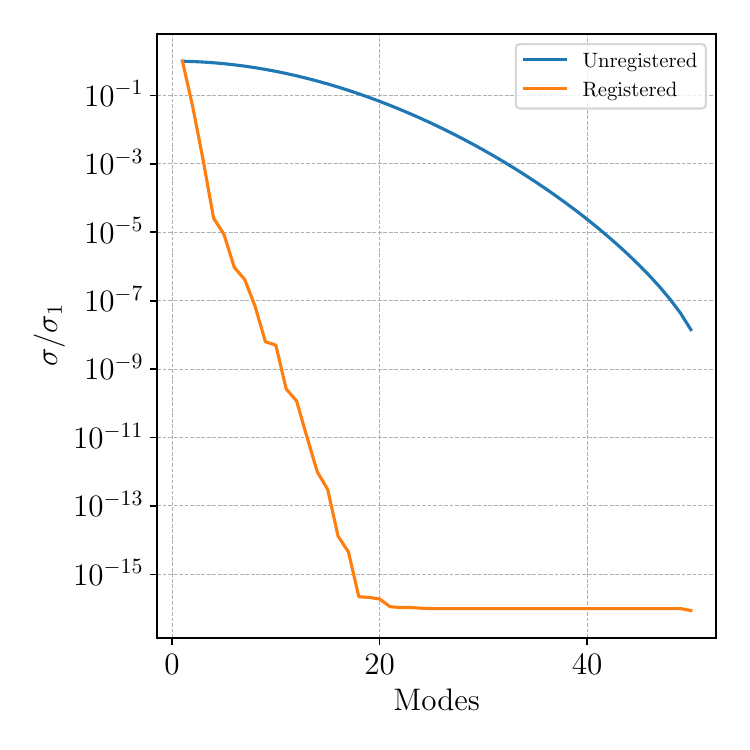}
        \caption{Singular values decay comparison}
        \label{fig:subfig1}
    \end{subfigure}
    \begin{subfigure}[b]{0.5\textwidth}
        \centering
        \includegraphics[width=\textwidth]{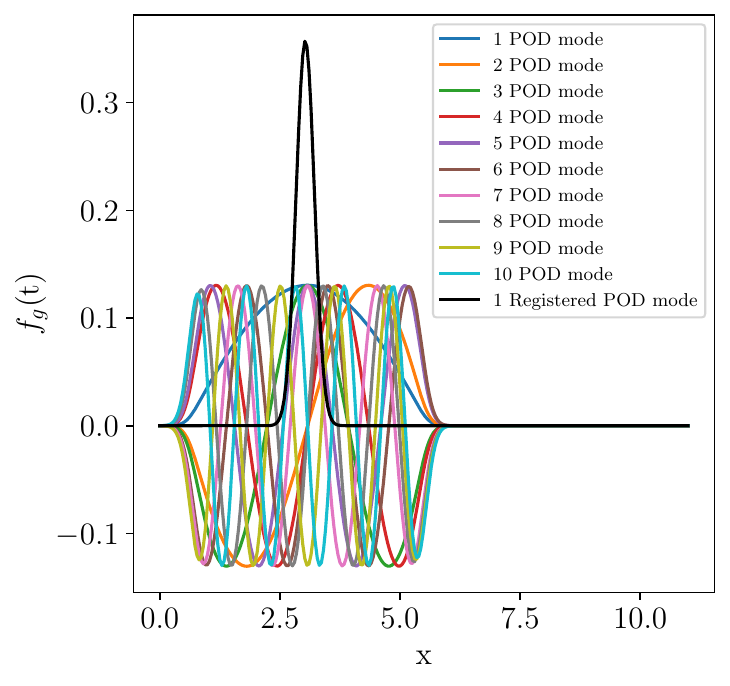}
        \caption{1D travelling wave registered POD mode}
        \label{fig:1D_registered_POD_mode}
    \end{subfigure}
    \caption{\textbf{Left:} Comparison of the singular values of registered and unregistered 1D travelling waves. \textbf{Right:} The first 10 POD modes obtained from the unregistered snapshots, alongside the 1st POD mode obtained via cross-correlation based snapshot registration.}
    \label{fig:1D_singular_values}
\end{figure}

In the offline phase, we select the reference snapshot corresponding to the parameter $t=3.03$. Using the procedure outlined in Algorithm \ref{Algorithm1}, we perform snapshot registration, aligning multiple snapshots to the reference by leveraging cross-correlation based registration, as shown in Figure \ref{fig:1D_registered_snapshots}. From Figure \ref{fig:1D_singular_values} it is evident that the registration significantly accelerates the KnW decay, leading us to seek a low-dimensional linear approximation subspace for efficient ROM construction. Figure \ref{fig:1D_registered_POD_mode} presents the first POD mode obtained by employing registration, which captures the structure of the Gaussian travelling waves, alongside the first 10 POD modes of the unregistered snapshots. In the case of unregistered snapshots, we have to consider numerous POD modes to construct ROM, which leads to inefficient ROMs. We construct ROM by considering the first POD mode (with 99.99 $\%$ of cumulative energy) obtained via registration, and we exploit the POD-R strategy discussed in section \ref{ROM}. Here, we utilize the radial basis function (RBF) interpolation for the regression task. \\

\begin{figure}[h!]
    \centering
    \includegraphics[scale=0.35]{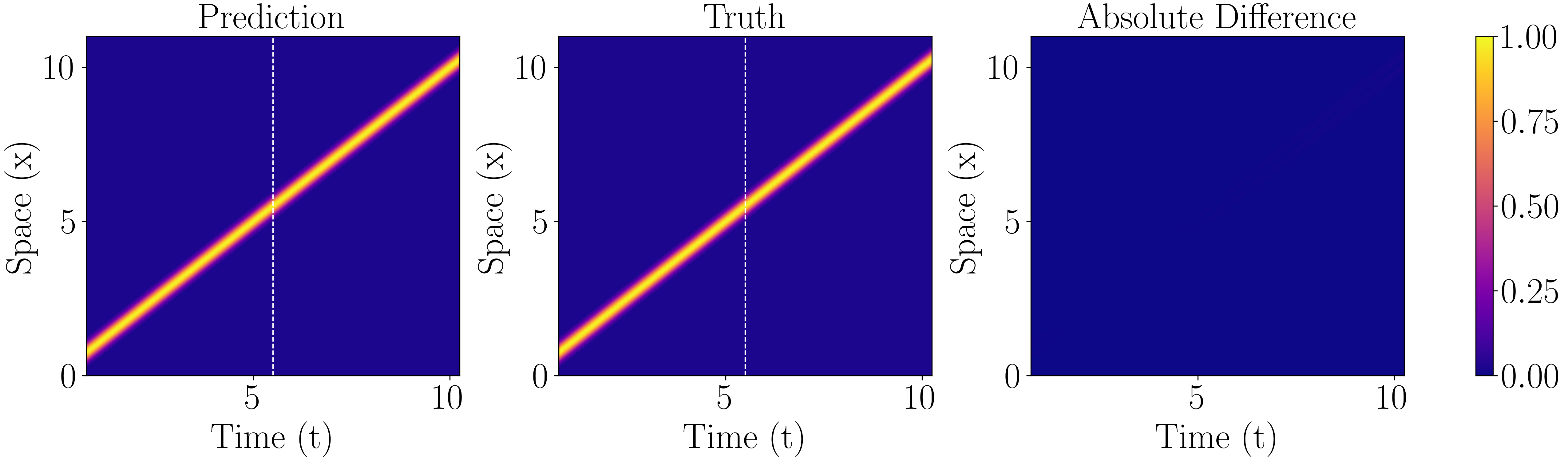}
    \caption{Comparison of the predictions of the resulting ROM of the 1D travelling waves with the original snapshots as ground truth, here the white line separates the first 50$\%$ dataset as the training set and the remaining as the test set and also shows the absolute difference between the prediction and the ground truth.}
    \label{fig:1Dpredcited_snapshots}
\end{figure}

\begin{figure}[h!]
    \centering
    \includegraphics[scale=0.5]{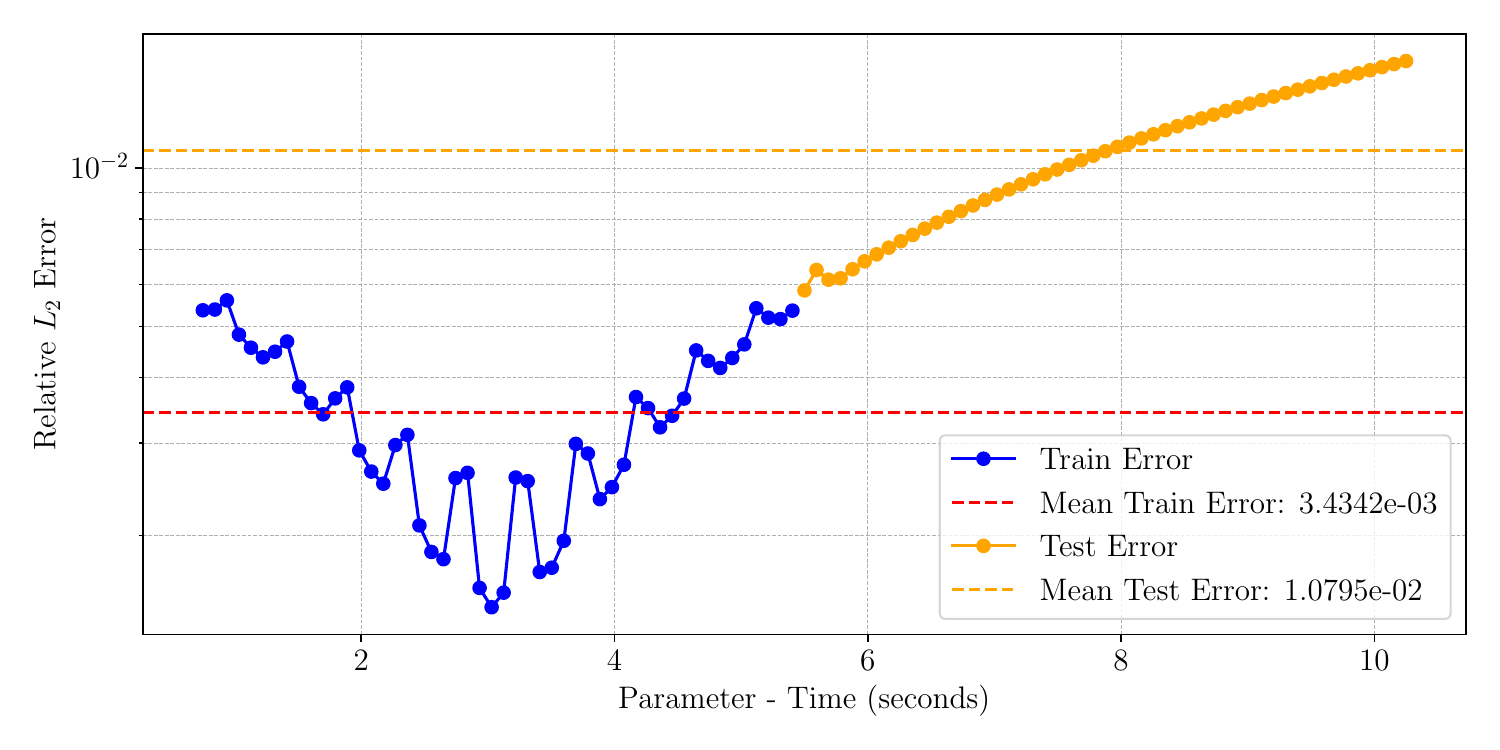}
    \caption{Relative $L_2$ predictions error by the ROM of the 1D travelling wave via cross-correlation based snapshots registration, for both training set and test set parameters.}
    \label{fig:1Dprediction_relative_error}
\end{figure}

After constructing the ROM, we incorporate shifts in the online phase by mapping parameter \( t_i \) to the optimal shift \( \Delta_i \) (denoted as \( \alpha \) in Algorithm \ref{Algorithm2}) using linear interpolation. To predict snapshots at new time instances, we employ the regression \( \alpha \) to transform the ROM-predicted snapshots (in the reference frame) back to the correct physical frame. Figure \ref{fig:1Dpredcited_snapshots} presents the predicted snapshots in the correct physical frame for both train and test set parameters, alongside the original snapshots (ground truth) and the absolute differences between predictions and ground truth. The relative prediction error, shown in Figure \ref{fig:1Dprediction_relative_error}, reveals that test errors are an order of magnitude larger than training errors, with the mean relative error highlighted. \\

\begin{figure}[h!]
    \includegraphics[scale=0.6]{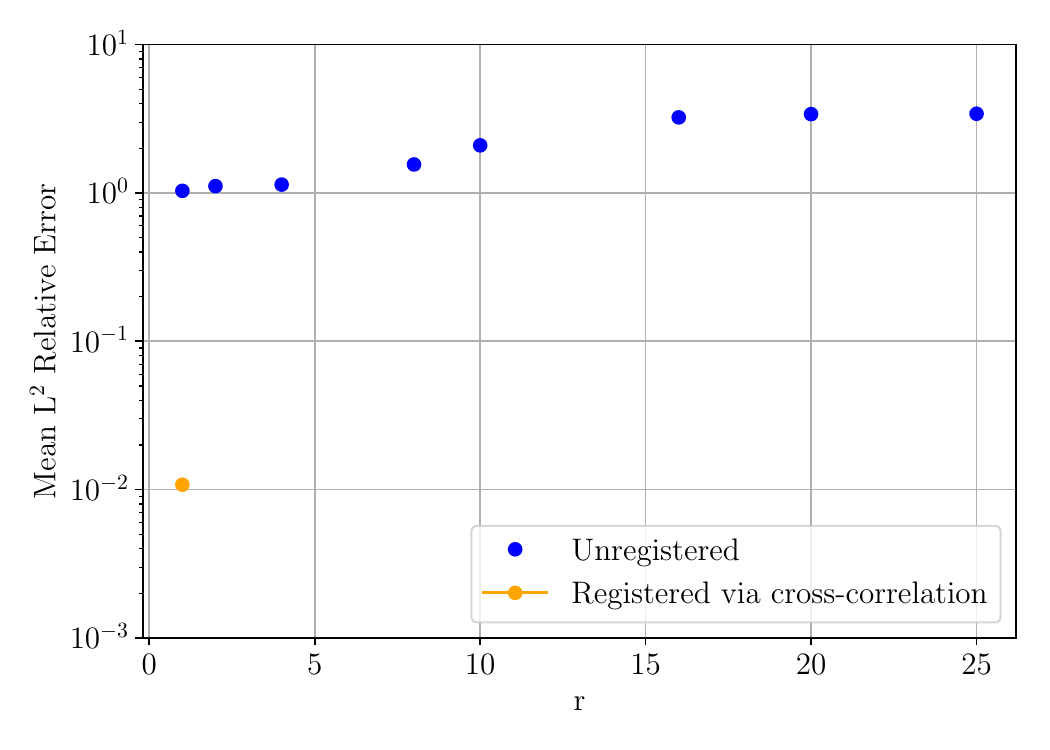}
    \centering
    \caption{Comparison of mean relative $L_2$ prediction error for 1D travelling wave test set parameters using POD-RBF-ROM, presented for the unregistered case (without registration) and the registered case (with registration via cross-correlation).}
    \label{fig:1D_mean_prediction_relative_error_test_set}
\end{figure}

In Figure \ref{fig:1D_mean_prediction_relative_error_test_set}, we present the mean relative $L_2$ prediction error for test set parameters, with varying POD ranks, for the construction of POD-RBF-ROM, both unregistered and registered snapshots via cross-correlation.  It is evident that in the unregistered case, even after considering 25 modes (with 99.99 $\%$ cumulative energy), it leads to inefficient and inaccurate ROMs even for simple linear transport problems, with errors two orders of magnitude larger than in the registered case. In contrast, for the registered case, the first mode (retaining 99.99$\%$ of the cumulative energy) results in an efficient and accurate ROM.

\subsection{2D isentropic convective vortex}

For the second experiment, we consider the 2D isentropic convective vortex, a popular higher-order method benchmark case \cite{spiegel2015survey, wang2013high}. To describe the flow physics of the inviscid, compressible flow of our benchmark case, we consider the 2D Unsteady Euler equations: \\ 

\begin{equation}
\frac{\partial}{\partial t}\left[\begin{array}{c}
\rho \\
\rho u \\
\rho v \\
\ E
\end{array}\right]+\frac{\partial}{\partial x}\left[\begin{array}{c}
\rho u \\
\rho u^2 + p \\
\rho u v \\
u(E+p)
\end{array}\right]+\frac{\partial}{\partial y}\left[\begin{array}{c}
\rho v \\
\rho u v \\
\rho v^2 + p\\
v(E+p)
\end{array}\right]=0
\label{eqn:Euler equation}
\end{equation}

where, $\rho$ is density, ($u,v$) are velocity x and y component, $p$ is pressure and $E$ is total energy. In this test case, we consider constant specific heat ratio, $\gamma = 1.4$, and gas constant, $R_{gas} = 287.15$ J/Kg.K. The thermodynamic closure is attained by the total energy equation (\ref{eqn: Energy equation}).

\begin{equation}
    E = \frac{p}{\gamma-1} + \frac{1}{2}\rho ( u ^{2} + v^{2}).
    \label{eqn: Energy equation}
\end{equation}
\begin{equation}
\left\{\begin{array}{l}
    \rho = \left[ 1 - \frac{\left(\gamma-1\right)b^2}{8\gamma\pi^2} e^{1-r^2} \right]^{\frac{1}{\gamma-1}}, \\
    p = \rho^\gamma, \\
    u = u_\infty - \frac{b}{2\pi} e^{\frac{1}{2}\left(1-r^2\right)} \left(y-y_c\right), \\
    v = v_\infty + \frac{b}{2\pi} e^{\frac{1}{2}\left(1-r^2\right)} \left(x-x_c\right). \\
\end{array}\right.
\label{Eqn: initial condition vortex}
\end{equation}

The freestream flow is initialized with the conditions \(\rho_\infty = 1\), \(u_\infty = 0.1\), \(v_\infty = 0\), and \(p_\infty = 1\). The vortex is introduced into the domain based on the analytical solution described in (\ref{Eqn: initial condition vortex}), with a vortex strength of \(b = 0.5\), a vortex radius \(r = \sqrt{(x - x_c)^2 + (y - y_c)^2} = 0.5\), and a vortex centre located at \((x_c, y_c) = (5, 10)\). There are several vortex problems reported in the literature \cite{ChiWang1997ENO-WENO}. In this study, we adopt the version implemented in the HyPar 1.0 package — finite difference hyperbolic-parabolic PDE solver on cartesian grids \cite{Hypar1} and we use this package to obtain the full-order solution of the equation \ref{eqn:Euler equation}. \\ 

To approximate the spatial derivatives in the governing equations, we employ the fifth-order hybrid-compact weighted essentially non-oscillatory (WENO) scheme, a high-order accurate numerical method designed for hyperbolic partial differential equations and convection-dominated problems \cite{zhang2016ENO-WENO}. For time integration, we consider the third-order strong stability-preserving Runge-Kutta (SSPRK3) method \cite{gottlieb2001SSPRK}. \\

\begin{figure}[h!]
    \centering
    \includegraphics[scale=0.55]{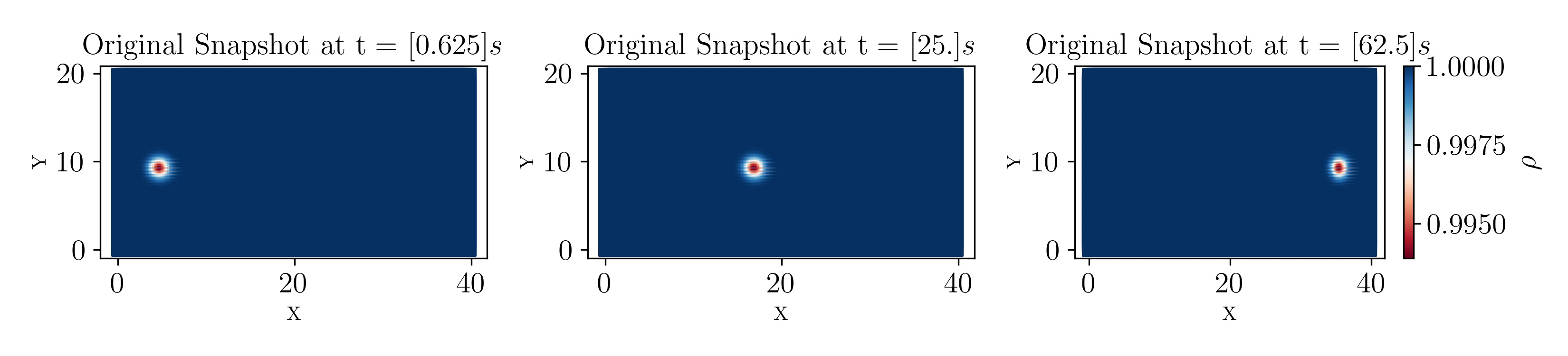}
    \caption{ Full-order solutions of the isentropic convective vortex at time instances, $ t \in \{0.625, 25, 62.5 \}$. }
    \label{fig:vortex_original_snapshots}
\end{figure}

We consider the rectangular computational domain \(\Omega = [0, 40] \times [0, 20]\), discretized into 28,800 cartesian grid points.  To introduce a vortex into the domain, we use the initial conditions given in (\ref{Eqn: initial condition vortex}). The simulation is performed over the time interval \(t \in [0, 62.5]\) using a time step of \(\Delta t = 0.00625\). Snapshots are recorded every 100 time-steps $\Delta t$, resulting in a snapshot matrix containing density field snapshots at \(N_T = 100\) time instances. The full-order solutions (snapshots) at three different time instances are shown in Figure \ref{fig:vortex_original_snapshots}. \\

\begin{figure}[h!]
    \centering
    \includegraphics[scale=0.45]{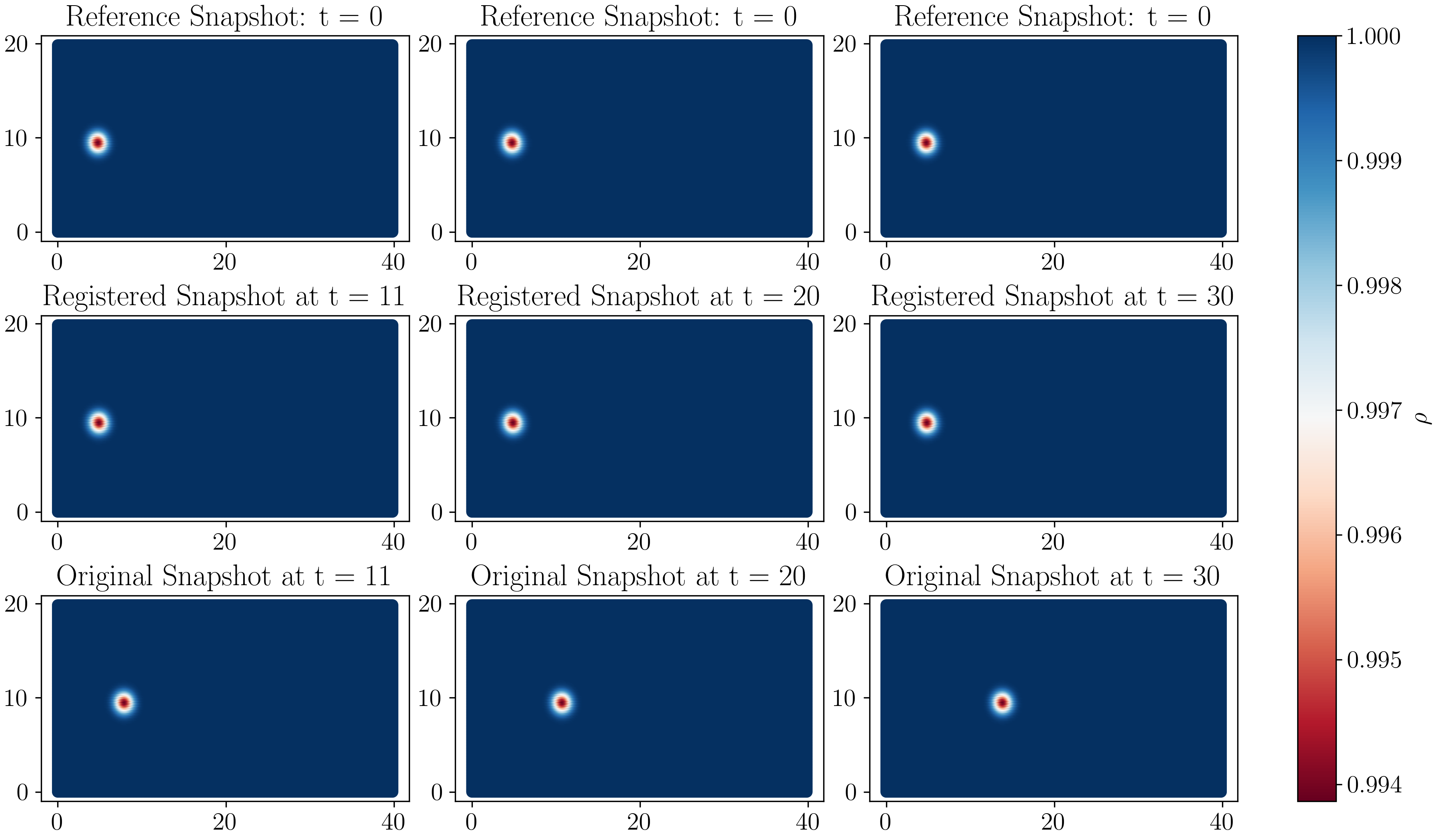}
    \caption{Comparison of the isentropic convective vortex registered snapshots with the reference and the corresponding original snapshot.}
    \label{fig:vortex_registered_snapshots}
\end{figure}

\begin{figure}[h!]
    \centering
    \begin{subfigure}[b]{0.4\textwidth}
        \includegraphics[width=\textwidth]{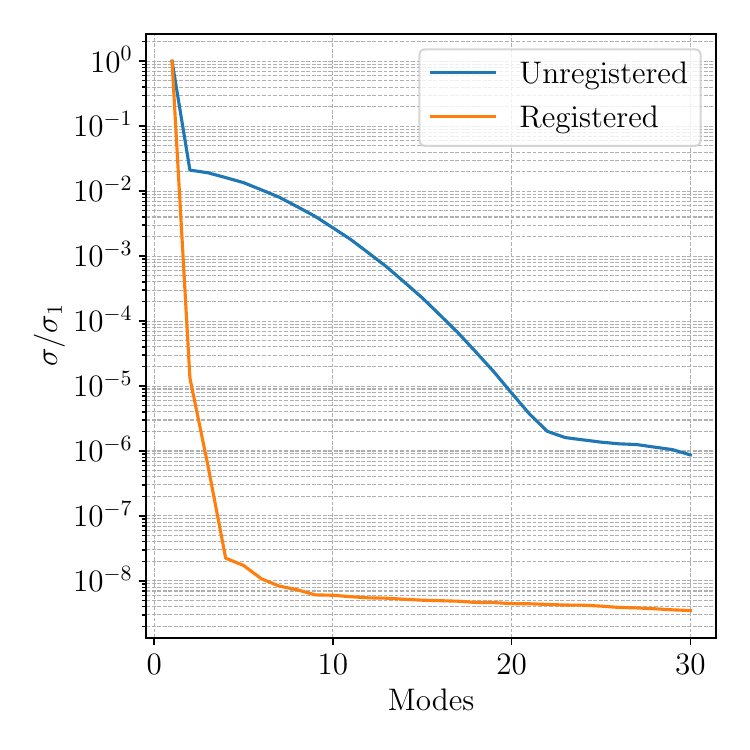}
        \caption{Singular values decay}
        \label{fig:vortex_singular_values_decay}
    \end{subfigure}
    \begin{subfigure}[b]{0.5\textwidth}
        \includegraphics[width=\textwidth]{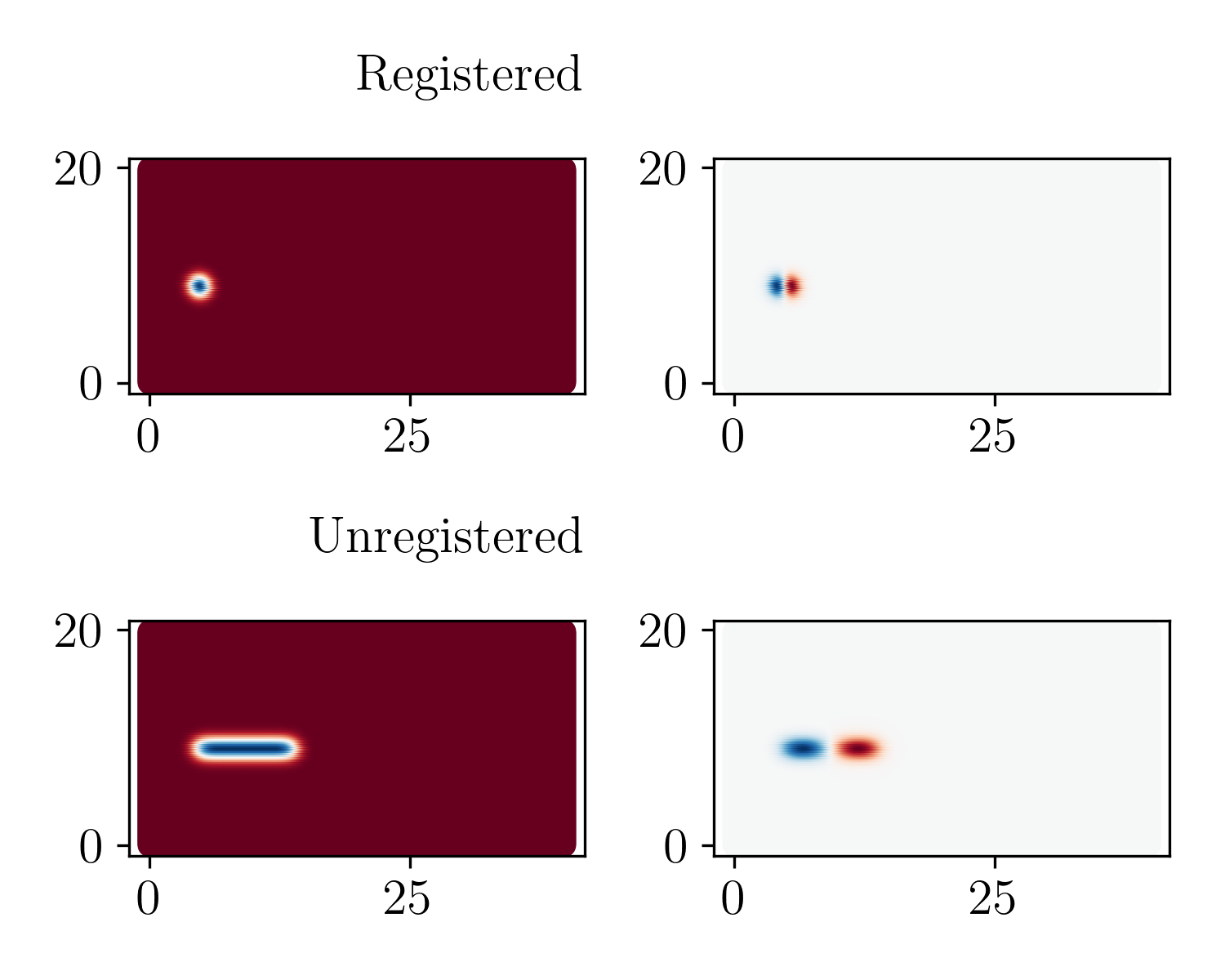}
        \caption{Comparison of POD modes}
        \label{fig:vortex_pod_modes}
    \end{subfigure}
    \caption{ \textbf{Left}: Comparison of singular values decay of registered snapshots matrix, labelled as Registered and of the unregistered snapshots matrix as Unregistered. \textbf{Right}: Showing the first two modes of the registered snapshots (first row - Registered) and showing the first two modes of the unregistered snapshots matrix (second row - Unregistered).}    
    \label{fig:vortex_sing_compare_POD_Mode}
\end{figure}

In this numerical experiment, we consider the first 30$\%$ dataset as a training set and the remaining 70$\%$ as the test set. We perform registration on the training set, Figure \ref{fig:vortex_registered_snapshots} presents the registered snapshots at time instances $ t \in \{0.625, 25, 62.5 \}$, comparing with the reference and the corresponding original snapshot. We perform reduction using POD on both the registered snapshots matrix and unregistered snapshots matrix, the singular values decay is shown in Figure \ref{fig:vortex_singular_values_decay} with sharp decay for registered snapshots. Furthermore, Figure \ref{fig:vortex_pod_modes} illustrates the POD modes obtained from the registered snapshots, which accurately capture the vortex structure (first row of Figure \ref{fig:vortex_pod_modes}). In contrast, the POD modes from the unregistered snapshots fail to capture the vortex structure effectively (second row of Figure \ref{fig:vortex_pod_modes}). \\

\begin{figure}[h!]
    \centering
    \includegraphics[scale=0.5]{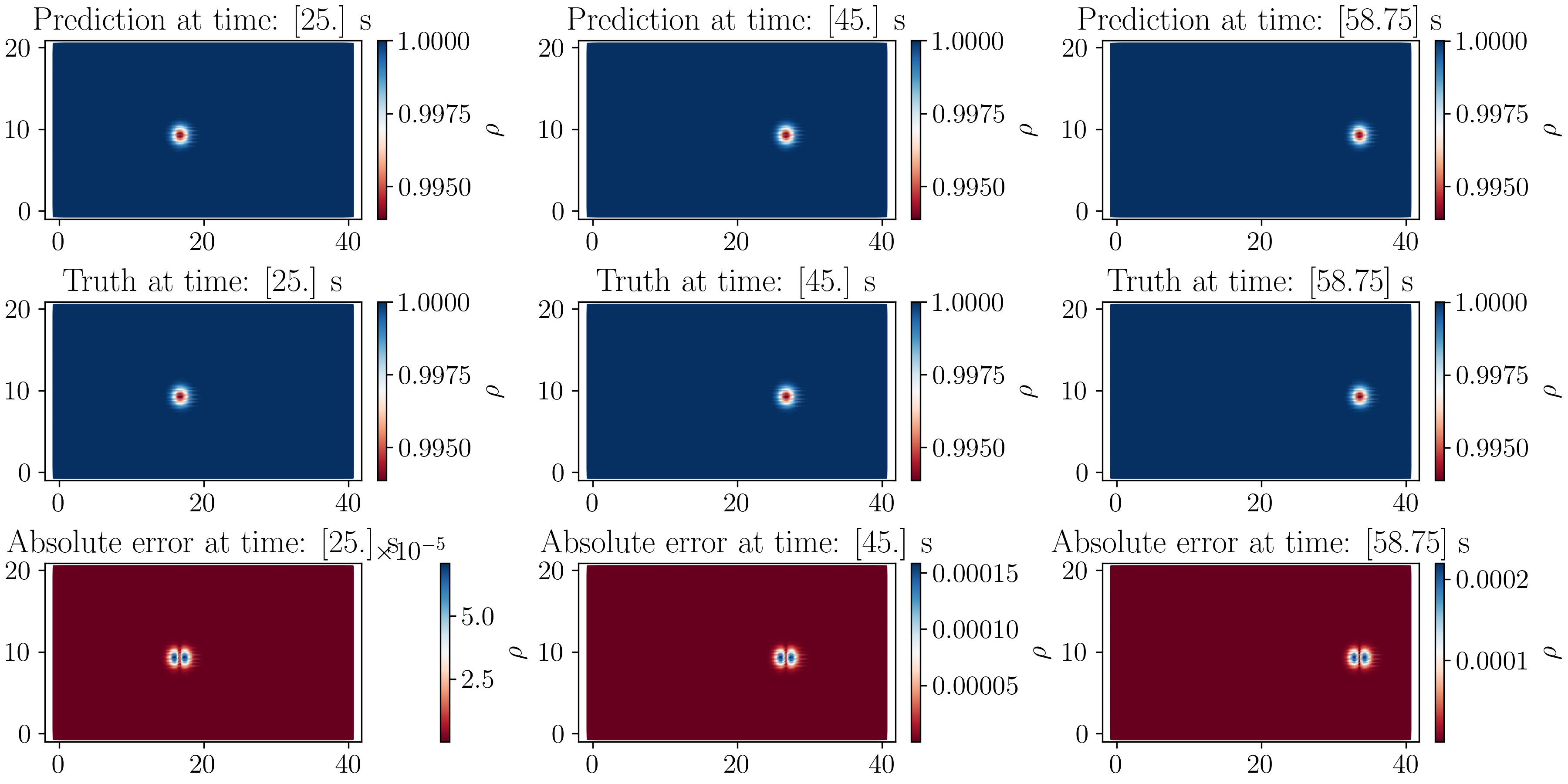}
    \caption{Predicted snapshots using POD-RBF-ROM with cross-correlation based registration for the isentropic convective vortex test set parameters at time instances $ t \in \{25, 45, 58.75 \}$s.}
    \label{fig:vortex_predcition_snapshots}
\end{figure}

\begin{figure}[h!]
    \centering
    \includegraphics[scale=0.45]{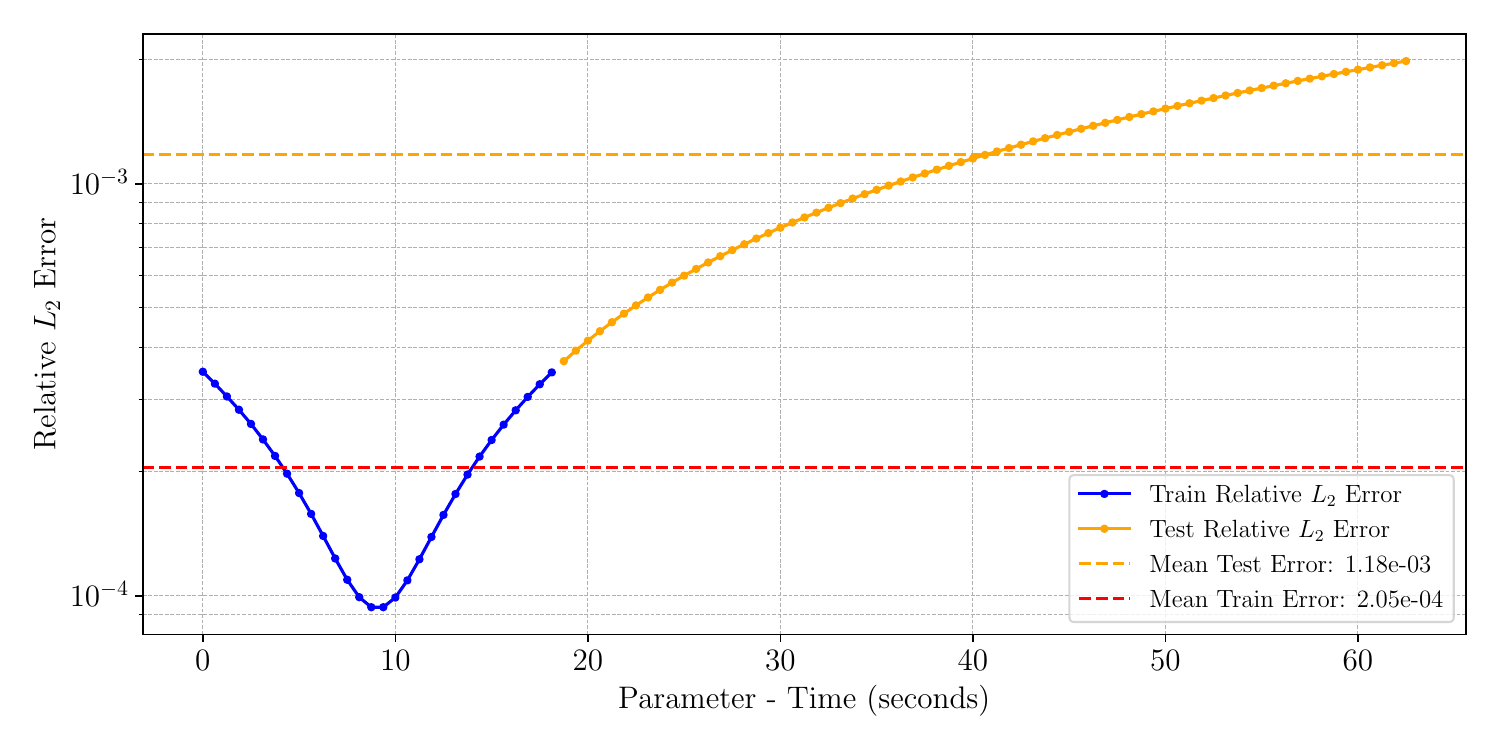}
    \caption{Relative $L_2$ predictions error of the isentropic convective vortex ROM via cross-correlation based registrations, for both training and test set parameters.}    
    \label{fig:Vortex_prediction_test_error}
\end{figure}

We utilize the first POD mode obtained via registration and as previously, we employ the POD-R strategy to construct ROM for this numerical experiment. To employ shifts in the online phase, we again use linear interpolations to build parameter \( t_i \) to optimal shift \( \Delta_i \) mappings. As specified in the previous case in the online phase, we utilize the predicted optimal shifts alongside the ROM prediction to revert the predicted snapshots to the correct physical frame. The predicted snapshots in the correct physical frame for three-time instances $ t \in \{25, 45, 58.75 \}$ from the test set are shown in Figure \ref{fig:vortex_predcition_snapshots}, with the corresponding original snapshot as ground truth and the absolute difference. The prediction error for both the training and test set parameters is presented in \ref{fig:Vortex_prediction_test_error}. The relative \( L_2 \) error of the predicted snapshots is plotted for each parameter, along with the mean relative \( L_2 \) error. \\

\begin{figure}[h!]
    \includegraphics[scale=0.5]{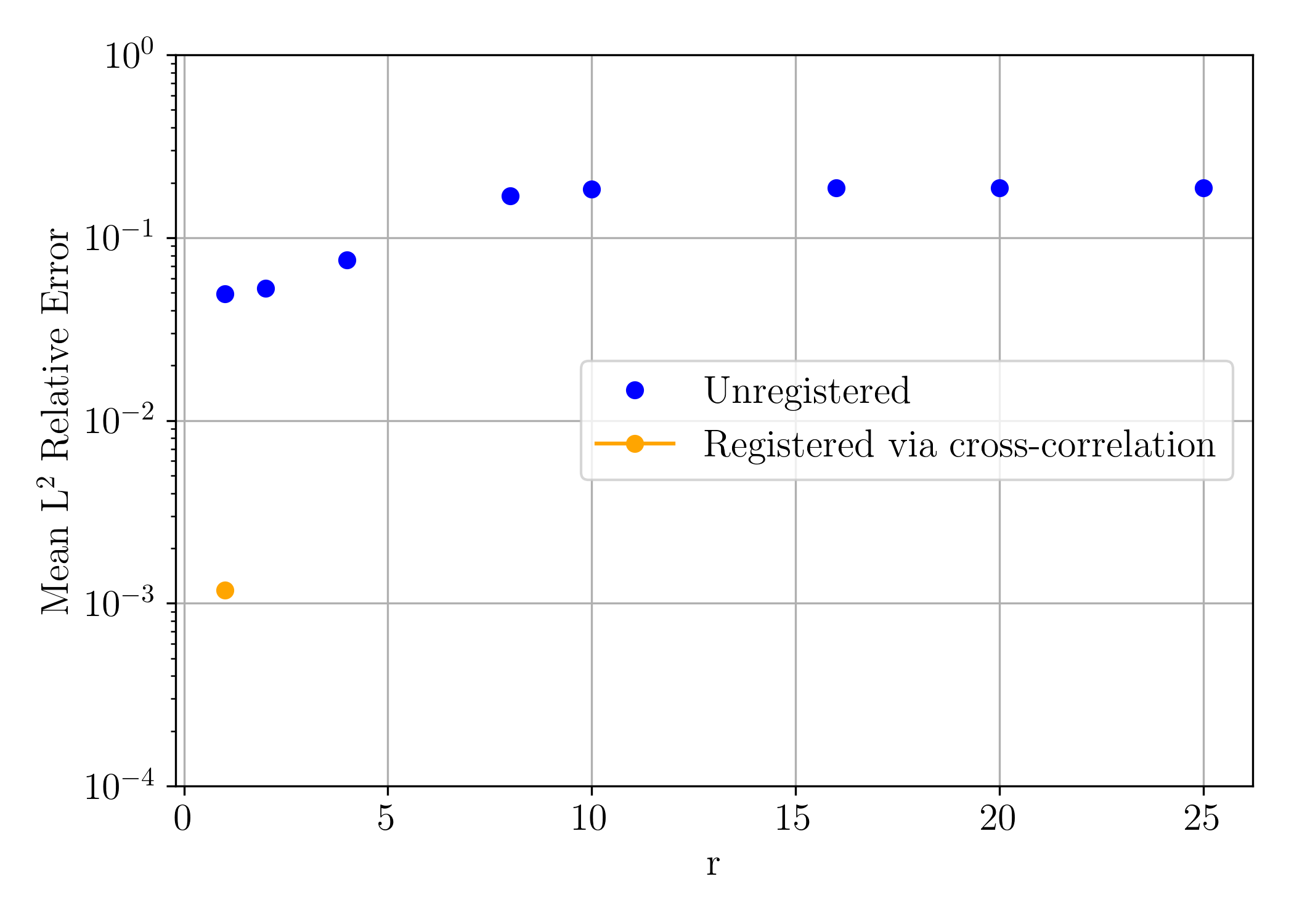}
    \centering
    \caption{Comparison of mean relative $L_2$ prediction error for 2D isentropic convective vortex test set parameters using POD-RBF-ROM, presented for the unregistered case (without registration) and the registered case (with registration via cross-correlation).}
    \label{fig:2D_vortex_prediction_relative_error}
\end{figure}

In Figure \ref{fig:2D_vortex_prediction_relative_error}, we illustrate the mean relative $L_2$ prediction error for test set parameters across different POD ranks in the construction of POD-RBF-ROM, using both unregistered and registered snapshots via cross-correlation. The results show that in the unregistered case, even with 25 modes, the ROM remains inefficient and inaccurate, yielding errors that are two orders of magnitude higher than in the registered case. Conversely, in the registered case, the first POD mode retaining 99.99 $\%$ of the cumulative energy is sufficient to yield efficient and accurate ROM.

\section{Conclusion}\label{sec:Conclusion}

In this work, we employ cross-correlation based snapshot registration for model reduction of transport-dominated problems. This cross-correlation based registration accelerates the KnW decay, facilitating the construction of a low-dimensional linear approximation subspace for efficient ROM development. We propose a complete framework comprising offline-online stages for the construction of ROMs using the cross-correlation based snapshots registration. We tested our proposed workflow on transport-dominated problems: 1D travelling waves, and 2D isentropic vortex. In both numerical experiments, the proposed approach results in efficient and accurate ROMs. However, the applicability of this approach is confined to certain transport-dominated problems where the structures of the transported features remain consistent. This approach is not suitable for addressing advection-diffusion problems. \\

Snapshot registration has seen significant advancements in medical imaging, largely driven by the progress of deep learning frameworks. Methods such as VoxelMorph \cite{Balakrishnan2019}, LapIRN \cite{Mok2020}, and FAIM (FAST Image Registration) \cite{Kuang2019} have demonstrated both speed and accuracy. These methods can be integrated within model reduction frameworks to address the challenges posed by slow KnW decay.

\section*{Acknowledgements}
We acknowledge the PhD grant supported by industrial partner Danieli \& C. S.p.A. and Programma Operativo Nazionale Ricerca e Innovazione 2014-2020, P.I. Gianluigi Rozza. GS acknowledges the financial support under the National Recovery and Resilience Plan (NRRP), Mission 4, Component 2, Investment 1.1, Call for tender No. 1409 published on 14.9.2022 by the Italian Ministry of University and Research (MUR), funded by the European Union – NextGenerationEU– Project Title ROMEU – CUP P2022FEZS3 - Grant Assignment Decree No. 1379 adopted on 01/09/2023 by the Italian Ministry of Ministry of University and Research (MUR) and acknowledges the financial support by the European Union (ERC, DANTE, GA-101115741). Views and opinions expressed are however, those of the author(s) only and do not necessarily reflect those of the European Union or the European Research Council Executive Agency. Neither the European Union nor the granting authority can be held responsible for them.

\section*{\large CRediT authorship contribution statement}
{Harshith Gowrachari}: Writing - original draft, Conceptualization, Data curation, Formal Analysis, Visualization, Methodology, Software. {Giovanni Stabile}: Writing – review $\&$ editing, Supervision. {Gianluigi Rozza}: Writing – review $\&$ editing, Funding acquisition, Project administration, Supervision. 

\section*{Data Statement}
The results presented can be reproduced by following the public repository at \newline
\href{https://github.com/harshith-gowrachari/MOR-CrossCorr-Registration}{https://github.com/harshith-gowrachari/MOR-CrossCorr-Registration}, which contains Jupyter notebooks for these numerical experiments. In this framework, we use the \texttt{scipy} library to compute cross-correlation and \texttt{EZyRB} package (\href{https://github.com/mathLab/EZyRB}{https://github.com/mathLab/EZyRB}) \cite{demo18ezyrb}, an open-source \texttt{Python} library for non-intrusive data-driven model order reduction of parameterized problems, developed and maintained at SISSA mathLab.
\bibliographystyle{abbrvnat} 
\bibliography{bib/biblio} 

\end{document}